\magnification=1100
\overfullrule=0mm

\input xypic
\input graphicx

\xyoption{all}
\xyoption{arc}
\let\s\scriptstyle

\def\Figure[#1]#2->#3;{  \vbox{    \halign{\parskip0pt   
\hfill ##   \hfill \cr \includegraphics[scale=#1]{#3.eps} \cr  \cr  \sl
 --  Figure  #2    --    \cr}  }
}

 \def\Figurea[#1]#2->#3;{  \vbox{    \halign{\parskip0pt   
\hfill ##   \hfill \cr \includegraphics[scale=#1]{#3.eps} \cr    \sl
 --  Figure  #2    --    \cr}  }
}

\def\point#1=#2;{\save   #1*{.}= "#2";  \restore}
\def\lettre#1,#2=#3;{\save   #2*+{#1}= "#3";  \restore}
\def \bpoint#1=#2;{\save   #1*-<2pt>{\bullet}= "#2";  \restore} 
\def \cpoint#1=#2;{\save   #1*-<2pt>{\circ}= "#2";  \restore} 
\def \vpoint#1=#2;{\save   #1*{}= "#2";  \restore} 
\def\aller#1-#2\par#3;{\save "#1"; "#2"**\crv{#3};\restore}
\def\daller#1-#2;{\save "#1"; "#2"**\dir{-};\restore}
\def\daller#1-#2;{\save "#1"; "#2"**\dir{-};\restore}
\def\trou#1=#2;{\save \POS?(#1)*{\hole}="#2";\restore}
\def\ttrou#1=#2;{\save \POS?(#1)*{o}="#2";\restore}
\def\otrou#1=#2;{\save \POS?(#1)*{\hole}*\frm{o}="#2";\restore}
\def \boite#1,#2=#3;{\save  #2*+{\s #1}*\frm{-}= "#3";  \restore}
\def \rond#1=#2;{\save  #1*{\hole}*\frm{o} = "#2";  \restore}
\def \unite#1=#2;{\save  #1*+<4.5pt>{\s \eta}*\frm{o} = "#2";  \restore}
\def \counite#1=#2;{\save  #1*+<4.5pt>{\s\varepsilon}*\frm{o} = "#2";  \restore}
\def \antip#1=#2;{\save  #1*+<5pt>{\s\sigma}*\frm{o} = "#2";  \restore}

\font\fiverom=cmr10 at 5pt
\font\sevenrom=cmr10 at 7pt
\font\eightrom=cmr10 at 8pt
\font\ninerom=cmr10 at 9pt

\font\twelvecmb=cmb10 at 12pt
\font\tencmb=cmb10 at 10pt
\font\tencmbi=cmbxti10 at 10pt
\font\ninecmb=cmb10 at 9pt

\font\tencmbxsl=cmbxsl10 at 10pt

\font\nineromtt=cmtt10 at 9pt

\font\bb=msbm10

\font\fivesymb=msam9 at 5pt

\font\tensymbo=msbm9

\font\sevensymbo=msbm9 at 7pt

\font\tensymbol=cmsy10

\font\sevensymbol=cmsy10 at 7pt

\font\calli=rsfs10 at 10pt

\def\HC{\hbox{$\cal H$}}

\def\HBB{\hbox{\bb H}}

\def\Hr{\hbox{\rm H}}

\def\calliC{\hbox{{\calli C}}}

\def\YD{\hbox{{\calli YD}}}

\def\HE{H{{\star}} E}

\def\ZBB{\hbox{\bb Z}}
\def\Fb{\hbox{\tencmb F}}
\def\Gb{\hbox{\tencmb G}}

\def\cross{\hbox{\tensymbol{\char 2}}}
\def\ZC{\hbox{$\cal Z$}}
 \def\sqtimes{\hbox{{$\cross$} \kern -1.26em
\raise 1.25ex \hbox { $\scriptstyle {\star}$}$\, $}}
\def\AC{\hbox{$\cal A$}}

\def\Hr{\hbox{\rm H}}

\def\Zr{\hbox{\rm Z}}

\def\Del {\hbox{$\Delta$}}
\def\id{\mathop{\hbox{\rm id}}\nolimits}
\def\idr{\hbox{\sevenrom id}}

\def\dm{ \null\hfill {\fivesymb {\char 3}}}
\def\lr {\hbox{$\ \longrightarrow\ $}}
\def\ot {\hbox{$\otimes$}}
\def\pa{\S\kern.15em}
\def\Dem{\noindent {\sl Proof:}$\ $}

\def\hfll#1#2{\smash{\mathop{\hbox to 10mm{\rightarrowfill}}
\limits^{\scriptstyle#1}_{\scriptstyle#2}}}
\def\cross{\hbox{\tensymbol{\char 2}}}
\def\lcross{\hbox{\sevensymbol{\char 2}}}

\def\sd{\hbox{\tensymbo{\char 111}}}
\def\lsd{\hbox{\sevensymbo{\char 111}}}

\newtoks\auteurcourant      \auteurcourant={\hfil}
\newtoks\titrecourant       \titrecourant={\hfil}

\newtoks\hautpagetitre      \hautpagetitre={\hfil}
\newtoks\baspagetitre       \baspagetitre={\hfil}

\newtoks\hautpagegauche   \newtoks\hautpagedroite

\hautpagegauche={\eightrom\rlap{\folio}\eightrom\hfil\the\auteurcourant\hfil}
\hautpagedroite={\eightrom\hfil\the\titrecourant\hfil\eightrom\llap{\folio}}

\newtoks\baspagegauche      \baspagegauche={\hfil}
\newtoks\baspagedroite      \baspagedroite={\hfil}

\newif\ifpagetitre          \pagetitretrue

\headline={\ifpagetitre\the\hautpagetitre
            \else\ifodd\pageno\the\hautpagedroite
             \else\the\hautpagegauche
              \fi\fi}

\footline={\ifpagetitre\the\baspagetitre\else
            \ifodd\pageno\the\baspagedroite
             \else\the\baspagegauche
              \fi\fi
               \global\pagetitrefalse}

\def\raggedbottom{\topskip 10pt plus 36pt\r@ggedbottomtrue}

\newcount\notenumber \notenumber=1
\def\note#1{\footnote{$^{{\the\notenumber}}$}{\eightrom {#1}}%
\global\advance\notenumber by 1}

\auteurcourant={Marc WAMBST}
\auteurcourant={Philippe NUSS -- Marc WAMBST}

\titrecourant={NON-ABELIAN HOPF COHOMOLOGY OF RADFORD PRODUCTS}

\noindent {\twelvecmb NON-ABELIAN HOPF COHOMOLOGY OF RADFORD PRODUCTS

}
\vskip 10pt

\bigskip

\noindent {PHILIPPE NUSS, MARC WAMBST}

\vskip 3pt
\noindent {\ninerom Institut de Recherche Math\'ematique Avanc\'ee,
Universit\'e de Strasbourg et CNRS, 7, rue Ren\'e-Descartes,
67084 Strasbourg Cedex, France.
e-mail: {\nineromtt nuss@math.u-strasbg.fr} and
{\nineromtt wambst@math.u-strasbg.fr}}

\vskip 15pt
\itemitem{}{\ninecmb Abstract.} {\ninerom
We study the non-abelian Hopf cohomology theory of Radford products
with coefficients in a comodule algebra.
We show that these sets can be expressed in terms of the non-abelian Hopf cohomology theory
of each factor of the Radford product. We write down an exact sequence relating these objects.
This allows to compute explicitly the non-abelian Hopf cohomology sets in 
large classes of examples.}
 \vskip 5pt
\itemitem{}\noindent {\ninecmb MSC 2000 Subject Classification.}
{\ninerom Primary: 18G50, 16W30 ; Secondary:   14A22, 16S35, 20J06, 55U10.}

\vskip 5pt
\itemitem{}\noindent {\ninecmb Key-words:} {\ninerom non-abelian cohomology, Radford products,
Hopf comodule algebra,  
cosimplicial non-abelian groups, Taft algebras.}

\vskip 20pt

\noindent {\tencmb I{\ninecmb NTRODUCTION}.}

\noindent The present article is devoted to the study of the non-abelian Hopf cohomology theory of Radford products.  The cohomology theory for Hopf algebras with coefficients in comodule algebras
was defined by the authors in  [4]. It generalizes the non-abelian cohomology theory for groups ([6],~[7]) and is closely related to
descent cohomology [3].
Recall that the Radford product  of two Hopf algebras (also called in 
the literature Radford biproducts
or bosonization) was introduced in [5] as a generalisation of the semi-direct product of groups (see also Majid's  interpretation in terms of braided 
categories~[2]).

The Radford product $\HE$, which is a Hopf algebra, is defined for a Hopf algebra $H$ and any 
algebra $E$ in the braided Yetter-Drinfeld category over  $H$ [11]. 
Our main result allows to  describe the cohomology sets of  $\HE$ in terms of those of $H$ and of  a braided version of the cohomology sets of $E$. This decomposition has a counterpart 
in the non-abelian cohomology theory of a semi-direct product of groups. We now detail the organisation of the article. 
\medskip

 In Section 0, we write down our conventions, for instance the diagrammatic language of Yetter-Drinfeld categories that we abundantly make use of. 
Section 1 is devoted to the Radford products. 
We briefly recall their definition and point out some properties of their comodule algebras. 
In~particular we show that a comodule algebra structure over a Radford product $\HE$ is equivalent to the datum of two compatible $H$- and $E$-comodule algebra structures (Proposition 1.1).

Reminders on the non-abelian Hopf cohomology theory $\HC^*(H, F)$ of a Hopf algebra $H$ with coefficients in a comodule algebra $F$ are given in Section 2. There we extend essentially in two ways the definition
introduced in [4]  in terms of pre-cosimplicial algebras.  Firstly we show that any functor~$\Gb$ defined on a suitable subcategory of the category of algebras with values in the latter 
gives rise to a new non-abelian Hopf cohomology theory
$\HC^*{\Gb}(H, F)$. Secondly we construct a theory $\HC^*_{\star}(E,F)$ by adapting the definitions 
when the Hopf algebra $E$ belongs to the braided Yetter-Drinfeld category of $H$.
The main result (Theorem~2.2) expresses  $\HC^*(\HE,F)$ in terms of $\HC^*(H, F)$ and 
$\HC^*_{\star}(E,F)$. We moreover show that the 1-cohomology set $\HC^1(\HE,F)$ fits into an exact sequence of pointed sets
$$ \xymatrix{{\cal H}^1{\bf C}_E(H, F)  \ar@<0ex>[r]^{  }
 & 
\ {\cal H}^1(\HE, F)  \ \  \ar@<0ex>[r]^{ } 
& \ {\cal H}^1_{\star}(E, F)^H,}$$
where the first map is injective (Theorem 2.4). The left-hand side term is the 1-cohomology set associated with  
the ``$E$-coinvariant" functor ${\bf C}_E = {\hbox {(--)}}^E$, and the right-hand side term is the $H$-coinvariant subset of 
${\cal H}^1_{\star}(E, F)$ in a sense we define
in the core of the text.
\goodbreak
Section 3 is devoted to the computation of some examples.
We begin with the case where the coefficient comodule algebra $F$ is equal to 
$E$. In Proposition 3.1 we show that  $\HC^0(\HE, E)$ is equal to the group of units of $k$,
and that $\HC^1(\HE, E)$ is in bijection with the set of groupoidal elements in $H$.
We specify this example when $\HE$ is a Taft algebra (Corollary 3.2). 
We also study the example when $H$ and $E$ are Hopf algebras of functions over groups (Proposition 3.3).
Finally, Section 4 is a foray into the world of groups.
We describe the non-abelian cohomology sets of a semi-direct product of two groups (Proposition~4.1) and obtain a similar decomposition as that stated in Theorem 2.2.

\bigskip
\noindent {\bf 0. Conventions and prerequisites.}

\smallskip

\noindent
Let $k$ be a fixed commutative and unital ring. By a {\sl $k$-module} we understand a symmetric {\hbox{$k$-bimodule.}}
The unadorned symbol $\ot$ between two {\hbox{$k$-modules}}  stands for $\ot _k$.
By {\sl (co-)algebra} we mean a (co-)unital (co-)associative {\hbox{$k$-(co-)}}al\-gebra.
By {\sl (co-)module} over a (co-)algebra $E$, we always understand a right $E$-(co-)module unless otherwise stated.
Let $M$ be a $k$-module. We identify in a systematic way $M \ot k$ with $M$. 
For any algebra $F$, we denote by $F^{\lcross}$ the group of invertible elements in $F$.
 
Let $H$ be a Hopf algebra  with multiplication $\mu_H$, unit map $\eta _H$,
comultiplication $\Delta _H$, counit map $\varepsilon _H$, and  antipode $\sigma_H$.
 To denote the coactions on elements, we use the Sweedler-Heyneman convention, that is for any Hopf algebra  $H$, the coproduct
$\Delta _H (h)$ of an element $h \in H$ is denoted by $h _1 \ot h_2$ (notice the  omission of the summation sign).
If $M$ is an $H$-comodule with coaction $\varrho _{M, H}$, for any $m \in M$ we set
$\varrho _{M, H}(m) = m_0 \ot m_1$.
From now on, we implicitely mean that $H$ is an $H$-module by $\mu _H$ and an $H$-comodule by $\Delta _H$.

In the sequel we shall use the diagrammatic language for morphisms (see for example [11]). The pictorial conventions are given in Figure 0.1.

\bigskip

$$ \Figure[ 1 ] 0.1->fig0-1; $$
 
\medskip
Notice that an element $m$ of a module $M$ is identified with the $k$-linear map from $k$ to $M$ sending $1$ to $m$. Hence that element is pictorially represented by a box without input and with $M$ as output.
Recall that, given two $k$-modules $M$ and $N$,  the {\sl flip}  is the map from $M\ot N$ to $N \ot M$ sending $m\ot n$ to
$n \ot m$, for any $m \in M$ and $n\in N$.

The tensor product $M\ot M'$ of two  $H$-modules $M$ and $M'$  becomes an $H$-module via the action given by
$(m\ot m'){\cdot}h = m{\cdot}h_1 \ \ot \ m'{\cdot}h_2$. Here $m$ belongs to $M$, $m'$ to $M'$, and $h$ to $H$.
Dually, the tensor product $N\ot N'$ of two  $H$-comodules $N$ and $N'$  is an $H$-comodule through the coaction $\varrho _{N{\otimes}N'\! , H}$ defined on $n\ot n' \in N\ot N'$ by
$\varrho _{N{\otimes}N'\! , H}(n\ot n')  = n_0 \ot n'_0 \ot n_1n'_1.$

Let $M$ still  be an $H$-module and $N$ be an $H$-comodule. In [11], Yetter defines a map
$\tau _{M, N} : M \ot N \lr N\ot M$ given, on any indecomposable tensor $m\ot n$ in $M \ot N$,  by 
$\tau _{M, N}(m\ot n) = n_0 \ot m{\cdot} n_1$. We shall use the following diagrammatic notation for $\tau _{M, N}$:

 \bigskip

   $$ \Figurea[ 1 ]0.2->fig0-2;  $$
   
   \bigskip
\noindent It is known that the map $\tau _{M, N}$ is a pre-braiding [11]. This means that the four equalities drawn down in Figure  0.3 hold, when  $M$ and $M'$ are $H$-modules,  $N$ and $N'$ are $H$-comodules,
 {\hbox{$\varphi : M \lr M'$}} is a morphism of $H$-modules, and $\psi : N \lr N'$  is a morphism of $H$-comodules.

\medskip
$$ \Figure[ 1 ]0.3->fig0-3;  $$

\medskip

Moreover, for any $H$-module $M$ and any $H$-comodule $N$, one easily shows the relations drawn in Figure~0.4.
We shall need them frequently in the sequel.
\medskip
$$ \Figure[ 1 ]0.4->fig0-4;  $$
\medskip

The Yetter-Drinfeld category $\YD_H^H$ of a  Hopf algebra $H$ is defined in the following way ([11]). Its objects are the
$k$-modules  $M$ which are both  $H$-modules and $H$-comodules subject to the compatibility law
$m_0 {\cdot} h_1 \ot m_1h_2  = (m {\cdot} h_2)_0 \ot h_1(m {\cdot} h_2)_1$,
 where $m \in M$ and $h \in H$. The latter relation is drawn in Figure 0.5.
 \medskip
$$\raise-55pt\Figure[ 1 ]0.5->fig0-5; $$
\medskip

\noindent {\bf Assumption :} {\sl From now on and to the end of the article, we suppose  that $H$ is a Hopf algebra such that the co-opposite algebra of $H$ is a Hopf algebra. }

\medskip
Under this assumption,  the pre-braiding $\tau _{-, -}$ is a braiding and $(\YD_H^H, \ot )$ is a braided strict monoidal category [11]. 
Remark that  the comodule structure of the tensor product $M \ot N$ of two $H$-comodules $M$ and $N$ is given by the formula
$  \varrho_{M\otimes N, H} = (\id _{M}\ot \tau _{H,N})\circ (\varrho_{M, H}\ot \id_{N}),$
summarized in Figure 0.6.

$$ \Figurea[ 1 ]0.6->fig0-6;  $$
\medskip

\noindent {\bf Lemma 0.1.} {\sl Let $H$ be a Hopf algebra. Let $E$ be a bialgebra in the category $\YD_H^H$ and let $N$ be an $H$-comodule. The map $\tau _{E, N} : E\ot N \lr N \ot E$ defined by Figure 0.2 is a morphism of $H$-comodules.}

\medskip
\Dem The proof is entirely contained in the pictorial computation of Figure 0.7.

\medskip
$$ \Figure[ 1 ]0.7->fig0-7;  $$
The first and the fourth equalities are consequences of Relation 0.4.b. The second one follows from the Yetter-Drinfeld compatibility (Figure 0.5). The third one is an obvious property of the flip, whereas the fifth is the convention adopted in Figure 0.2.
\dm

\medskip

\setbox1\hbox{\includegraphics[scale=1]{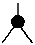}}
\setbox2\hbox{\includegraphics[scale=1]{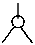}}

For any bialgebra $E$ in the category $\YD_H^H$, we define an algebra $F$ to be a  {\sl  Radford $E$-comodule algebra} if $F$ is both an $H$-comodule (pictorially  \raise-2pt\box1 \ ) and an $E$-comodule (pictorially \raise-2pt\box2 \ ), such
that the three conditions drawn in Figure 0.8 are satisfied.

\medskip
 
$$ \Figure[ 1 ]0.8->fig0-8;  $$
\medskip

\noindent In particular, any $E$-comodule algebra $F$ in $\YD_H^H$ -- for example $E$ itself -- is a Radford $E$-comodule algebra. Indeed, the two first relations are consequences of the $E$-comodule structure of $F$. The third one is the transcription of the fact that $\varrho _{F,E}$ is a morphism of $H$-comodules. We choose the name  {\sl Radford comodule algebra}
in view of the results obtained in Proposition 1.1 of the next section.



\bigskip
\smallskip\smallskip

\noindent {\bf 1. Radford products}

\smallskip

\smallskip

\noindent {\tencmbxsl 1.1. Some definitions.}
Recall that, for any   bialgebra $E$ in $\YD_H^H$, one may construct the Radford product 
$\HE$ ([5], [1]). It is a bialgebra and contains
$H$ as a sub-bialgebra and $E$ as a subalgebra.
As a module,  $\HE$ is equal to $H \ot E$. The product, coproduct, unit, and counit are respectively given on $h{\star} x$ and $h'{\star} x'$ in $\HE$ by
$$\eqalign{(h{\star} x)(h'{\star} x') =  hh'_1\star (x\cdot h'_2)x'  & \qquad \Delta_{\HE}(h{\star} x) =  \bigl(h_1{\star} (x_1)_0\bigr) \otimes \bigl(h_2(x_1)_1{\star} x_2\bigr)\cr
1_{\HE} = 1_{H}{\star }1_{E}& \qquad \varepsilon_{\HE}(h{\star }x) =   \varepsilon _H(h)  \varepsilon _E (x).\cr}$$
\goodbreak
The diagrammatic transcription of the two first formulas is drawn in Figure 1.1.
 \medskip
 
$$ \Figurea[ 1 ]1.1->fig1-1;  $$
\medskip
 The canonical projections and inclusions are summarised in
the following diagram:

$$ \xymatrix{ H\ar@<1ex>[rr]^{\idr_H \otimes \eta _E  }  & &  \HE 
   \ar@<1ex>[ll]^{\idr_H \otimes \varepsilon _E }
\ar@<1ex>[rr]^{\varepsilon _H \otimes \idr _E  }
 & & E \ar@<1ex>[ll]^{\eta _H \otimes \idr _E  }}$$
 \medskip
\noindent  The two left maps are morphisms of bialgebras, whereas   in general $\eta _H \otimes \id _E$ is a only a map of algebras and $\varepsilon _H \otimes \id _E$ is only a map of  coalgebras.

\medskip

We also recall that, if $E$ is a Hopf algebra in $\YD_H^H$, one may turn $\HE$ into a Hopf algebra by using the map given by Figure 1.2 
as an antipode [5].

$$ \Figurea[ 1 ]1.2->fig1-1a;  $$

When $E$ is trivial both as an $H$-module and as an $H$-comodule, then  $\HE$ is simply the tensor product $H \ot E$ of the two Hopf
algebras $H$ and $E$.

\bigskip
\noindent {\tencmbxsl 1.2. Comodule algebras over Radford products.}

\smallskip
\noindent {\bf Proposition 1.1.} {\sl  Let $H$ be a Hopf algebra. Let $E$ be a bialgebra in $\YD_H^H$ and $F$ be an algebra. The two following statements are equivalent:

\itemitem{1)} $F$ is an $\HE$-comodule algebra~;

\itemitem{2)} $F$ is an $H$-comodule algebra and $F$ is a Radford $E$-comodule algebra.}

\bigskip

\setbox1\hbox{\includegraphics[scale=1]{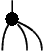}}
\setbox2\hbox{\includegraphics[scale=0.8]{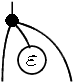}}
\setbox3\hbox{\includegraphics[scale=0.8]{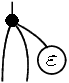}}
\Dem 
First suppose that $F$ is an $\HE$-comodule algebra. Let $\varrho _{F, H{\star} E} : F \lr 
F \ot (\HE)$ be the coaction of $H {\star} E$ on $F$
(pictorially \ \raise -3pt\box1 \ ). Define $\varrho _{F, H} : F \lr F \ot H $ and  $\varrho _{F, E} : F \lr F \ot E $ by 
$$\varrho _{F, H} = (\id _F \ot \id _H \ot \varepsilon _E)\circ \varrho _{F, H{\star} E}  \qquad {\hbox {\rm and}} \qquad
\varrho _{F, E} = (\id _F \ot \varepsilon _H \ot \id _E)\circ \varrho _{F, H{\star} E} $$

\noindent (pictorially \ \raise-6pt \box3 \ and \  \raise-6pt\box2 \ ).
\smallskip

The fact that $F$ is an $H{\star} E$-comodule is represented in Figure 1.3:

\bigskip
 
   $$ \Figure[ 1 ]1.3->fig1-2;  $$

\noindent We apply the counit maps $\varepsilon _H$, respectively $\varepsilon _E$, at the two corresponding outputs of the equality in Figure 1.3. Using the defining axiom
of the counit map and  the relations 0.3.c and 0.3.d of Figure~0.3, one gets the four relations of Figure 1.4.

\medskip
 
$$ \Figurea[ 1 ]1.4->fig1-3;  $$
\medskip
\noindent The two upper relations show that $\varrho _{F, H}$, respectively $\varrho _{F, E}$, endow $F$ with an $H$-, respectively {\hbox{$E$-, comodule}} structure
(in both cases, the axiom $(\id _F \ot \varepsilon _-)\circ \varrho _{F, -} = \id _F$ is obviously satisfied).
The two lower equalities imply the compatibility relation 0.8.c of Figure 0.8. Moreover, 
Relation~0.8.b is implied by
$\varrho _{F, -} \circ \eta _F = \eta _F \otimes \eta _-$, which is obvious.

Consider now Figure 1.5, which depicts that $\varrho _{F, H {\star} E}$ is a morphism of algebras.

 \medskip

$$ \Figure[ 1 ]1.5->fig1-4;  $$

\noindent Apply the counit map $\varepsilon _E$ to Figure 1.5. Using  the fact that $\varepsilon _E$ is a morphism of algebras and verifies Relation 0.3.c, one gets
 the following equality, which proves that $\varrho _{F,H}$ is a morphism of algebras:

 \medskip

$$ \Figure[ 1 ]1.6->fig1-5;  $$

\noindent Similarly, apply the counit map $\varepsilon _H$  to Figure 1.5. Using  the fact that $\varepsilon _H$ is a morphism of algebras and verifies Relations 0.4.c and 1.4.c, one obtains
 Relation 0.8.a for $\varrho_{F, E}$.
\medskip
$$ \Figure[ 1 ]1.7->fig1-6;  $$

    \medskip

\noindent This achieves the proof of 1 $\Longrightarrow$ 2.

Conversely, take  an $H$-comodule algebra $F$ (with coaction $\varrho _{F, H}$) which is also a
Radford  {\hbox{$E$-comodule algebra}}
(via the coaction $\varrho _{F, E}$).  We show that the map $\varrho _{F, H {\star} E} : F \lr F \ot (H {\star} E)$ defined by 
 $$\varrho _{F, H {\star} E} = (\varrho _{F, H} \ot \id _E)\circ \varrho _{F, E}$$
 endows $F$ with a structure of $\HE$-comodule algebra.
 
We obviously have the identity $(\id _F \ot \varepsilon _{H{\star}E})\circ \varrho _{F, H{\star}E} = \id _F$. The coassociativity
of $\varrho _{F, H{\star}E}$ results from the equalities 1.8.a of Figure 1.8 below. It is an immediate consequence of the compatibility relation 0.8.c of Figure 0.8
and of the definition of the comultiplication of $H{\star}E$. The fact that $F$ is an $H{\star}E$-comodule algebra follows from the obvious relation
$\varrho _{F, H{\star}E}\circ \eta _F = \eta _F \otimes \eta _{H{\star}E}$ and from the equalities 1.8.b. There one uses the axiom 0.8.a of Figure 0.8, the $H$-comodule algebra structure of $F$, and Relation 0.4.b of Figure 0.4. 

\medskip
$$ \Figure[ 1 ]1.8->fig1-7;  $$
\dm

\bigskip

\noindent {\tencmbxsl 1.3. Radford extensions of scalars.} 

\bigskip

\noindent {\bf Proposition 1.2.} {\sl Let  $H$ be a Hopf algebra. Let $E$ be a bialgebra in $\YD_H^H$ and  $F$
be a Radford $E$-comodule algebra. The module $F \ot E$ is a Radford $E$-comodule algebra for the multiplication given by Relation 1.9.a, the 
 unit $\eta _F \ot \eta _E$, the $H$-coaction map depicted in Relation 1.9.b, and the $E$-comodule map ${\id}_F \ot \Delta _E$.  }

\medskip
$$ \Figure[ 1 ]1.9->fig1-8;  $$

\medskip 

\Dem Using the diagrammatic language,  the Reader may check all the axioms. \dm

\medskip

The Radford $E$-comodule algebra of Proposition 1.2 will be denoted by $F {\star} E$. One iterates the construction in order to build $(F {\star} E){\star} E$. The latter Radford $E$-comodule algebra will be denoted
for simplicity by $F {\star} E {\star} E$. 

\medskip

\noindent {\sl Remark 1.3.~:} The coaction map $\varrho _{F,E} : F \lr F {\star} E$ is a morphism of $H$-comodules by Relation~0.8.c in Figure 0.8. 

\goodbreak
\bigskip
\smallskip

\noindent {\bf 2. Radford products and non-abelian Hopf cohomology theory.}

\smallskip
\noindent In this section, we recall and generalize  the construction of
the non-abelian Hopf cohomology theory developped in [4]. We study 
the case when the Hopf algebra is a Radford product.

\bigskip
\noindent {\tencmbxsl 2.1. Remainders on non-abelian Hopf cohomology theory.} 
By a {\sl pre-cosimplicial object in a category $\calliC$}, we understand a diagram $\AC$ of the
type $$\AC  =  \xymatrix{ A^0  \ar@<1.3ex>[r]^{d^0}
\ar@<-1.3ex>[r]^{d^1}  & \ A^1
\ar@<2.3ex>[r]^{d^0}
\ar@<0ex>[r]^{d^1}
\ar@<-2.3ex>[r]^{d^2}& \ A^2},$$ where $A^i$ is an object in $\calliC$ and $d^j$ are morphisms in $\calliC$ such that the
pre-cosimplicial relations $d^id^j = d^jd^{i-1}$  hold for $i > j$.

\medskip

\noindent Suppose now that $\calliC$ is the category of groups. The {\sl non-abelian $0$-cohomology
group} $\HBB^0(\AC)$ of a pre-cosimplicial group $\AC$ is the equalizer
of the pair $(d^0, d^1)$:
$$\HBB^0(\AC) = \{ x \in A^0 \ \vert \ d^1(x) = d^0(x) \}.$$
 The {\sl non-abelian $1$-cohomology pointed set} $\HBB^1(\AC)$ of $\AC$  is the right quotient
$$\HBB^1(\AC) =  A^0 \big\backslash \ZBB^1(\AC),$$
where the set $\ZBB^1(\AC)$ of {\sl $1$-cocycles}
is the subset of $A^1$ defined by
$$\ZBB^1(\AC) =\{ X \in A^1 \ \vert \ d^2(X)d^0(X) = d^1(X) \}.$$
The group $A^0$
acts on  the right on $A^1$  by $$X \leftharpoonup x = (d^1x^{-1}) X (d^0x),$$
where $X \in A^1$ and $x \in  A^0$. This action restricts to  $\ZBB^1(\AC)$.
Two $1$-cocycles $X$ and $X'$
are  {\sl cohomologous} if they belong to the same orbit
under this action. The quotient set $\HBB^1(\AC) =   A^0 \big\backslash \ZBB^1(\AC) $
is pointed with distinguished point the class of the neutral element of $  A^1$.
\medskip

\goodbreak
 Let $H$ be a Hopf algebra and $F$ be an $H$-comodule algebra
with multiplication $\mu_F$ and coaction
$\varrho _{F, H}$. In the spirit of [9], we defined in [4] two maps $d^i: F \lr F \ot H$
($i= 0,1$)
and three maps $d^i:  F \ot H \lr  F \ot H\ot H$ ($i=0,1,2$) 
by the formulae
$$\eqalign{d^0 (x)  &  = \varrho _{F, H} (x), \hskip47.55pt d^1 (x)  = x\ot 1,  \cr
 d^0(X)  & = (\varrho _{F, H} \ot \id _H) (X), \quad
 d^1(X)    = (\id _F \ot  \Delta _H)(X),  \quad   \ d^2(X)    = X \ot 1, \cr}$$
where $x \in F$ and $X \in F \ot H$. The maps $d^i$ are easily seen to be morphisms of algebras. The  diagram ${\cal C}(H, F)$ given by
$$  \xymatrix{ F  \ar@<1.3ex>[r]^{d^0 \ \ \ }
\ar@<-1.3ex>[r]^{d^1 \ \ \ }  & \ F \ot H
\ar@<2.3ex>[r]^{d^0\ \ \  }
\ar@<0ex>[r]^{d^1\ \ \ }
\ar@<-2.3ex>[r]^{d^2\ \ \ }& \ F \ot H \ot H}$$
is a pre-cosimplicial object in the category of algebras. In the sequel the diagram ${\cal C}(H, F)$ lives in fact in a certain subcategory $\cal B$
of the category of algebras, for instance a category of comodule algebras.

Let $\Gb$ be a functor from the category $\cal B$ to the category of algebras.
The diagram ${\cal C}(H, F)$ gives rise to a  pre-cosimplicial diagram ${{\cal C}^{\lcross}}{\Gb}(H, F)$ in the category of groups  by setting:


$${{\cal C}^{\lcross}}{\Gb}(H, F) \ = \  \xymatrix{ {\Gb} (F)^{\lcross}  \ar@<1.3ex>[r]^{d^0 \ \ \ }
\ar@<-1.3ex>[r]^{d^1 \ \ \ }  & \  {\Gb}(F \ot H)^{\lcross}
\ar@<2.3ex>[r]^{d^0 \ \ \ }
\ar@<0ex>[r]^{d^1 \ \ \ }
\ar@<-2.3ex>[r]^{d^2\ \ \ }& \ { \Gb}(F \ot H \ot H)^{\lcross}}$$
(we still denote by  $d^i$ the restrictions of the maps ${\Gb}(d^i)$ to invertible elements).
The {\sl non-abelian Hopf cohomology theory  $\HC^*{\Gb}(H, F)$ of $H$
with coefficients in $F$ and with respect to ${\Gb}$}
is  the non-abelian cohomology of the diagram  ${\cal C}^{\lcross}{\Gb}(H, F) $. In other words, we have
$$ \HC^0\Gb(H, F)  = \HBB^0 \bigl({\cal C}^{\lcross}\Gb(H, F) \bigr)  \quad {\hbox{\rm  and}}
\quad
\HC^1\Gb(H, F) = \HBB^1 \bigl({\cal C}^{\lcross}\Gb(H, F) \bigr).$$  
\noindent We denote  $\ZC^1{\Gb}(H, F)$ the set
$\ZBB^1\bigl({\cal C}^{\lcross}\Gb(H, F) \bigr)$ and call it
the set of Hopf $1$-cocycles of $H$ with coefficients in $F$.

\smallskip
In particular, the non-abelian cohomology theory $\HC^*(H, F)$ defined in [4] is the theory with respect to the  identity functor.  \medskip

We have now to add a technical lemma, which gives an alternative description of the $1$-cocyle set  $\ZC^1(H, F)$. This result is the counterpart in the cohomology theory of Proposition 2.1 of [4]
which was applied to the restricted cohomology theory studied in [3].

\medskip

\noindent {\bf Lemma 2.0:} {\sl For any $H$-comodule algebra, one has $$\ZC^1(H, F) =\{ X \in F\ot H \quad \vert \quad (\id_F \ot \varepsilon_H)(X) = 1 \quad {\sl and} \quad  d^2(X)d^0(X) = d^1(X) \}.$$ }

\Dem Pick an element $X \in F \ot H$ satisfying the cocycle condition. We have to prove that the invertibility of $X$ is equivalent to the equality  $(\id_F \ot \varepsilon_H)(X) = 1$. 
In the following, we set $x = (\id_F \ot \varepsilon_H)(X)$.

First suppose that $X$ is invertible in $F \ot H$.
 Applying the map $\id_F \ot \varepsilon_H \ot \id _H$ to the
cocycle relation, one gets the equality $(x\ot 1)X = X$ in $F \ot H$, hence $x\ot 1 =1 \ot 1$. By use of $\id_F \ot \varepsilon_H$, one deduces $x = 1$.

Conversely, suppose that $x$ is equal to $1$. Set $Y = ((\id _F \ot \mu _H)\circ (\varrho_F \ot \sigma _H))(X)$. Apply $(\id _F \ot \mu _H)\circ (\id_F \ot \id _H \ot \sigma _H)$ to the cocycle
condition $d^2(X)d^0(X) = d^1(X) $. The first picture of Figure 2.1 computes the left-hand side of the obtained equality whereas the second one calculates
the right-hand side.

\medskip
$$ \Figure[ 1 ]2.1->fig2-0;  $$
So, one gets $XY = x \ot 1 = 1 \ot 1$. 
In order to prove the left invertibility of $X$, apply the map
$(\id _F \ot \mu _H)\circ (\varrho_{F, H} \ot \mu _H ) \circ (\id _F \ot \sigma _H \ot \id _H)$ to the cocycle
condition. This leads to the Hopf yoga of Figure 2.2.
\medskip
$$ \Figure[ 0.9]2.2->fig2-00;  $$
This gives $YX = x \ot 1 = 1 \ot 1$.
\dm

\bigskip
\noindent {\tencmbxsl 2.2. Cohomology theory with coefficients in Radford comodule algebras.}
Let  $H$ be a Hopf algebra. Let $E$  be a Hopf algebra in the category $\YD_H^H$ and let $F$  be a   Radford $E$-comodule algebra with
$\varrho _{F, E} $ as $E$-comodule map.
We paraphraze the above constructions in order to construct ``braided" cohomology sets $\HC^*_{\star}(E, F)$.

One checks that the maps $d^*$ defined by the formulae
$$\eqalign{d^0 (x)  &  = \varrho _{F, E} (x), \hskip47.55pt d^1 (x)  = x\ot 1,  \cr
 d^0(X)  & = (\varrho _{F, E} \ot \id _E) (X), \quad
 d^1(X)    = (\id _F \ot  \Delta _E)(X),  \quad   \ d^2(X)    = X \ot 1 \cr}$$
(where $x \in F$ and $X \in F \ot E$) induce the following pre-cosimplicial diagram ${\cal C}_{\star}(E, F) $ in the category of algebras:
$$ \xymatrix{ F \  \ar@<1.3ex>[r]^{d^0 \ \ \ }
\ar@<-1.3ex>[r]^{d^1 \ \ \ }  &  \ F{\star} E \
\ar@<2.3ex>[r]^{d^0 \ \ \ }
\ar@<0ex>[r]^{d^1 \ \ \ }
\ar@<-2.3ex>[r]^{d^2\ \ \ }& \ F{\star}E {\star}E.}$$
\medskip
For any endofunctor $\Gb$ of the category of algebras, one defines in the same way as before the non-abelian Hopf cohomology theory  $\HC_{\star}^*\Gb(E, F)$ with respect to $\Gb$ to be the non-abelian cohomology objects associated with the
 pre-cosimplicial diagram of groups ${\cal C}_{\star}^{\lcross}{\Gb}(E, F)$. 
 The set of $1$-cocycles are denoted by $\ZC_{\star}^1\Gb(E, F)$.

 \smallskip
The non-abelian Hopf cohomology objects $\HC_{\star}^*(E, F)$
are  by definition $\HC_{\star}^*\Gb(E, F)$, with $\Gb$ being the identity functor.
Notice that {\sl mutatis mutandis} Lemma 2.0 is still valid for $\ZC_{\star}^1(E, F)$.
In~the case when $H$ is equal to the ground ring $k$, then $E$ is a Hopf algebra and $F$ is an $E$-comodule algebra in the usual sense and
 one has the equality $\HC_{\star}^*(E, F) = \HC^*(E, F)$.

\goodbreak 
\bigskip
\noindent {\tencmbxsl 2.3. The decomposition Theorem.} We are now able to state our main result, that is how to decompose, under the same hypotheses as before, the sets $\HC^*(H{\star}E, F)$ in terms of $\HC^*(H, F)$ and $\HC_{\star}^*(E, F)$.
The difficulty occures of course in the $1$-cohomology level. 

At this point we have to introduce
the intermediate set $\ZC^1(H, F) \sqtimes \ZC_{\star}^1(E, F)$ defined
as follows. An element $[U^H,U^E]$ of $\ZC^1(H, F) \sqtimes \ZC_{\star}^1(E, F)$ is a couple $(U^H,U^E)$ in $\ZC^1(H, F) \cross \ZC_{\star}^1(E, F)$ such that $U^H$ and $U^E$ satisfy the compatibility relation drawn down in Figure 2.3.
\medskip
 $$ \Figurea[ 1 ]2.3->fig2-1;  $$
\medskip
\noindent 
Observe that $\ZC^1(H, F) \sqtimes \ZC_{\star}^1(E, F)$ is a pointed set with $[1\ot 1, 1\ot 1]$ as distinguished point.

\medskip

\noindent {\bf Lemma 2.1:} {\sl Let $H$ be a Hopf algebra. Let $E$  be a Hopf algebra in the category $\YD_H^H$ and let $F$  be a Radford $E$-comodule algebra.
Then the group $F^{\lcross}$ acts on the right on  $\ZC^1(H, F) \sqtimes \ZC_{\star}^1(E, F)$ by
$$[U^H,U^E] \leftharpoonup x = [U^H \leftharpoonup x,U^E\leftharpoonup x],$$
where $[U^H,U^E]$ belongs to $\ZC^1(H, F) \sqtimes \ZC_{\star}^1(E, F)$ and $x$ is an element of $F^{\lcross}$.}

\bigskip

\Dem Take $[U^H,U^E]$ in $\ZC^1(H, F) \sqtimes \ZC_{\star}^1(E, F)$ and $x$ in $F^{\lcross}$. One has to show that the elements  $U^H \leftharpoonup x$ and $U^E\leftharpoonup x$ are compatible in the sense of Figure 2.3.
By definition, $U^H \leftharpoonup x$ and $U^E\leftharpoonup x$ are respectively given by 

$$\eqalign{U^H \leftharpoonup x & = (d^1x^{-1})U^H(d^0x) = (x^{-1}\ot 1)U^H(\varrho_{F,H}(x)) \ \ {\rm and} \cr
U^E \leftharpoonup x & = (d^1x^{-1})U^E(d^0x) = (x^{-1}\ot 1)U^E(\varrho_{F,E}(x)). \cr}$$
The product is  the multiplication in $F\ot H$ in the first line, and the multiplication in $F {\star} E$ in the second line.
These formulae may be rewritten as shawn in Figure 2.4.

$$ \Figurea[ 1 ]2.4->fig2-2;  $$\medskip

 Computing the left-hand side (respectively the right-hand side) of Relation 2.3 for the couple $(U^H \leftharpoonup x,U^E\leftharpoonup x)$, we obtain Figure 2.5  (respectively Figure 2.6)
 
 $$ \Figure[ 0.9 ]2.5->fig2-3;  $$
 \medskip
$$ \Figure[ 0.9 ]2.6->fig2-4;  $$\medskip

Let us explain  Figure 2.5. The first equality is a consequence of Figure 2.4. The second one comes from Relation 0.3.c applied to the $H$-module morphism $\mu _E$ and from the associativity of $\mu _F$ and of $\mu _E$.
The third equality is implied by Relation 0.3.d applied to the $H$-comodule morphism~$\mu _F$ and by the associativity of $\mu _F$. Equality (4) is deduced from Relation 0.8.a. Finally the fifth equality is the
simplification of the product $xx^{-1}$, which appears by associativity of $\mu _F$.

We detail now Figure 2.6.  The first equality is again a consequence of Figure 2.4. Equality (2) uses   Relation 0.3.d for $\mu _F$, the associativity of $\mu _F$ and also the simplification $xx^{-1}=1$. The third equality
comes from 0.3.d applied to the $H$-comodule morphism $\mu _E$ and from the associativity of~$\mu _F$. Fourth equality is the sum of 0.3.d for $\mu _F$ and the associativity of $\mu _F$.
Equality (5) derives from Lemma 0.1 and from Relation 0.3.d. The sixth equality  is the compatibility relation between $U^H$ and $U^E$ (Figure 2.3). Equality (7) translates the associativity of $\mu_F$
and 0.3.d for $\varrho _{F,E}$. The eighth equality holds for the same reasons as the third one.
Finally, by Relation 0.8.a, the two last pictures of Figure 2.5
and of Figure 2.6 are equal. \dm

\bigskip

\noindent {\bf Theorem 2.2:} {\sl Let $H$ be a Hopf algebra. Let $E$  be a Hopf algebra in the category $\YD_H^H$ and let $F$  be a Radford $E$-comodule algebra.
Then one has the equality
$$\HC^0(H{\star} E, F) = \HC^0(H, F) \cap \HC_{\star}^0(E, F),$$
and an isomorphism of pointed sets
$$\HC^1(H{\star} E, F) \cong F^{\lcross} \big\backslash (\ZC^1(H, F) \sqtimes \ZC_{\star}^1(E, F)),$$
where the group $F^{\lcross}$ acts on $\ZC^1(H, F) \sqtimes \ZC_{\star}^1(E, F)$ by
$$[U^H,U^E] \leftharpoonup x = [U^H \leftharpoonup x,U^E\leftharpoonup x].$$

}

\Dem {\sl $0$-cohomology level}. Recall that, after Proposition 1.1, the relationships between the three coaction maps on $F$ are given by

{\centerline{\vbox{\hsize300pt 
\halign{$#$ \hfill & $#$ & $#$ & \qquad # \hfill\cr
\varrho _{F, H{\star} E} &=& (\varrho _{F, H}\ot \id _E) \circ \varrho _{F,  E} = (\id _F \ot \tau _{E,H}) \circ (\varrho _{F, E}\ot \id _H) \circ \varrho _{F,  H} & (1)\cr
\varrho _{F, H}  &=& (\id _F \ot \id _H \ot \varepsilon _E) {\circ} \varrho _{F, H{\star} E} & (2)\cr
\varrho _{F, E}  &=& (\id _F \ot \varepsilon _H \ot \id _E) {\circ} \varrho _{F, H{\star} E} & (3)\cr}}}}

\smallskip
\noindent Let $x$ be an element of $\HC^0(H{\star} E, F)$, that is an invertible element of $F$ such that $\varrho _{F, H{\star} E} (x) = x\ot 1 \ot 1.$
The equalities (2) and (3) applied to $x$ give $\varrho _{F, H} (x) = x\ot 1 $ and $\varrho _{F, E} (x) = x\ot 1$. This means that $x$ belongs  simultaneously to
$\HC^0(H, F)$ and to $\HC_{\star}^0(E, F)$. Conversely, let $x$ be an element of $\HC^0(H, F) \cap \HC_{\star}^0(E, F)$. The first equality in (1) gives 
$$\varrho _{F, H{\star} E}(x) = ((\varrho _{F, H}\ot \id _E)\circ \varrho _{F,  E}) (x) = (\varrho _{F, H} \ot \id _E)(x \ot 1) = x \ot 1 \ot 1.$$
Hence $x$ is a $0$-cocycle of $H{\star} E$ with coefficients in $F$.

\medskip

\noindent {\sl $1$-cocycle level}. We first show that $\ZC^1(H{\star} E, F)$ is isomorphic to $\ZC^1(H, F) \sqtimes \ZC_{\star}^1(E, F)$. Let $X$ be an element of $\ZC^1(H{\star} E, F)$.
We define $X^H \in F \ot H$ and $X^E \in F \ot E$ by
$$X^H = (\id _F \ot \id_H \ot \varepsilon _E) (X) \quad {\rm and} \quad X^E = (\id _F \ot \varepsilon _H \ot \id_E) (X).$$
The element $X$ being invertible, so is $X^H$, since  $\id_H \otimes \varepsilon _E$ is a morphism of algebras. 
The invertibility of $X^E$ is a consequence of Lemma 2.0. Indeed, one has $(\id_F \ot \varepsilon_E) (X^E) = (\id_F \ot \varepsilon_{\HE}) (X) = 1$.

By definition, $X$ satisfies the relation $d^2(X)d^0(X) = d^1(X)$, which pictorially looks like Relation  2.7.a, or  after a slight simplification, like Relation  2.7.b.

$$ \Figurea[ 0.9 ]2.7->fig2-5;  $$\medskip

The four manners to apply the counities to the $H$- or $E$-outputs of 2.7.b. give rise to the four relations collected in Figure 2.8. 
 
$$ \Figurea[ 0.9 ]2.8->fig2-6;  $$\medskip

Relation 2.8.a  (respectively 2.8.d)
says nothing else that $X^H$ (respectively $X^E$)  belongs to $\ZC^1(H, F)$ (respectively $\ZC_{\star}^1(E, F)$). Relations 2.8.b and 2.8.c imply that $X^H$ and $X^E$ are 
compatible in the sense of Figure 2.3.

\medskip
Conversely, pick $[X^H,X^E]$ in $\ZC^1(H, F) \sqtimes \ZC_{\star}^1(E, F)$. Set 
$$X = \bigl((\mu_F \ot \id _H \ot \id_E) \circ (\id _F \ot \tau _{H,E} \ot \id _E)\bigr)(X^H \ot X^E),$$
as depicted in Figure 2.9.

 $$ \Figure[ 1 ]2.9->fig2-7;  $$\medskip
 
Figure 2.10 below shows that $X$ verifies the cocycle condition. The first equality is a consequence of Relation 0.8.a. The second and the third one come from the associativity of the products 
$\mu _H$ and $\mu _E$. Equality (4) is implied by Relation 0.3.c applied to $\mu _E$ and by Relation 0.3.d applied to $\mu _F$. The fifth equality is the Compatibility relation 2.3. The next one is a consequence of Relation 0.3.d successively applied to $\varrho_{F, E}$, to $\tau_{E,F}$ and to $\mu _E$.
The seventh equality is still a consequence of Relation 0.3.d applied to $\mu _F$. Equality (8) is the cocycle conditions for the two elements $X^H$ and $X^E$. Finally, the last one is
Relation 0.3.c applied to $\Delta _H$.

$$ \Figure[ 0.8 ]2.10->fig2-8;  $$\medskip

 In order to prove the invertibility of $X$, apply the map $\id_F \ot \varepsilon_{\HE}$ to $X$. By Relation 0.3.c and  Lemma~2.0, the result is equal to $1$,
hence, again by Lemma~2.0, $X$ belongs to $\ZC^1(H{\star}E, F)$. 
\medskip
\noindent {\sl $1$-cohomology level}. 
We leave to the Reader the proof of the two equalities 
$$(X \leftharpoonup x)^H = (X^H \leftharpoonup x) \quad {\rm and} \quad (X \leftharpoonup x)^E = (X^E \leftharpoonup x),$$
where $X \in \ZC^1(H{\star} E, F)$ and $x \in F^{\lcross}$. The first one is straightforward and the second one 
is a consequence of the Relations 0.4.c and 1.4.c.
\dm

\bigskip
\noindent {\tencmbxsl 2.4. An exact sequence of $1$-cohomology sets.} 
In this paragraph, we  construct under the previous hypotheses
an exact sequence of pointed sets associated with
the diagram 
$$ \xymatrix{ H\ar@<0ex>[rr]^{\idr_H \otimes \eta _E  }  & &  \HE 
\ar@<0ex>[rr]^{\varepsilon _H \otimes \idr _E  }
 & & E .}$$
 
 To this end, we introduce  the two following functors. 
 First let ${\bf C}_E$ be the endofunctor of the category of $E$-comodule algebras defined on a $E$-comodule algebra $F$
 by ${\bf C}_E(F) = F^E$. Here $F^E$ stands for the algebra of $E$-coinvariant elements of $F$. Secondly, let ${\bf T}_H$ be the
  endofunctor of the  category of 
 algebras given on an algebra $F$ by ${\bf T}_H(F) = F \ot H$.
 
Notice that for any non-negative integer $i$, the module $F\ot H ^{\otimes i}$ can be endowed with the structure of $E$-comodule algebra through the formula
 $\varrho _{F\otimes H^{\otimes i}\!, E}(x \ot {\underline {h}}) =  x_0 \ot \underline {h}\ot x_1$. Here $\underline {h}$ belongs to $H^{\otimes i}$,
 the element $x$ belongs to $F$, and we set $\varrho_{F, E}(x) = x_0 \ot x_1 $. With respect to this \hbox{$E$-comodule} structure, 
 ${\cal C}(H,F)$ may easily be seen to be a pre-cosimplicial diagram in the category of \hbox{$E$-comodule} algebras. One is therefore allowed to form the cohomology sets
 ${\cal H}^*{\bf C}_E(H, F)$. 
  \medskip
\noindent {\bf Lemma 2.3:}  {\sl Let $H$ be a Hopf algebra. Let $E$  be a Hopf algebra in the category $\YD_H^H$ and let $F$  be a Radford $E$-comodule algebra. When $H$ is flat as a $k$-module, there is an isomorphism }
$${\cal H}^*{\bf C}_E(H, F) \cong {\cal H}^*(H, F^E).$$
\noindent \Dem By flatness,  $(F\ot H ^{\otimes i})^E$ is isomorphic to  $F^E\ot H^{\otimes i}$, for $i = 0, 1, 2$. This implies an isomor\-phism
${\cal C}^*{\bf C}_E(H, F) \cong {\cal C}^*(H, F^E)$
on the co-simplicial  level, hence on the cohomology sets.\dm

\bigskip
As explained in 1.3, the algebras $F$, $F{\star} E$ and $F{\star} E{\star} E$ are $H$-comodule algebras.
Consider the two morphisms of  cosimplicial algebras
$ \xymatrix{ {\cal C}_{\star}(E, F) \ \  \ar@<0.5ex>[r]^{\eta \ \ }
\ar@<-0.5ex>[r]_{\varrho \ \  }  & \  \ {\cal C}_{\star}{\bf T}_H(E, F)}$
summarised in the following diagram:
$$ \xymatrix{ {\cal C}_{\star}(E, F)    \ar@<0.3ex>[dd]^{\varrho \ \ \ \ }
\ar@<-1.3ex>[dd]_{\eta }
 & = &
F \ \  \ar@<1.3ex>[rr]^{d^0 \ \ \ }\ar@<0.3ex>[dd]^{\varrho_{F, H} \ \ \ }
\ar@<-1.3ex>[dd]_{{\rm id} _F \otimes \eta _H  }
\ar@<-1.3ex>[rr]^{d^1 \ \ \ }  && \  \ F {\star} E \ 
\ar@<2.4ex>[rr]^{d^0 \ \ \ }
\ar@<0ex>[rr]^{d^1 \ \ \ }
\ar@<-2.4ex>[rr]^{d^2\ \ \ } \ar@<0.3ex>[dd]^{\varrho_{F{\star} E, H} \ \ \ }
\ar@<-1.3ex>[dd]_{{\rm id} _{F{\star}E} \otimes \eta _H  } && \ \  F {\star} E {\star} E 
\ar@<0.3ex>[dd]^{\varrho_{F{\star} E {\star} E, H} \ \ \ }
\ar@<-1.3ex>[dd]_{{\rm id} _{F{\star}E{\star}E} \otimes \eta _H  }
\\ \\ 
{\cal C}_{\star}{\bf T}_H(E, F) & = & F \ot H \ \  \ar@<1.3ex>[rr]^{d^0 \ \ \ }
\ar@<-1.3ex>[rr]^{d^1 \ \ \ }  &&  \ \  (F {\star} E) \ot H\
\ar@<2.4ex>[rr]^{d^0 \ \ \ }
\ar@<0ex>[rr]^{d^1 \ \ \ }
\ar@<-2.4ex>[rr]^{d^2\ \ \ } && \ \  (F {\star} E {\star} E) \ot H.}$$
They give rise to the two morphisms $ \xymatrix{ {\cal H}_{\star}^*(E, F) \ \  \ar@<0.5ex>[r]^{\eta \ \ \ \ }
\ar@<-0.5ex>[r]_{\varrho \ \ \ \ }  & \  \ {\cal H}_{\star}^*{\bf T}_H(E, F)}$ on the cohomology level.
Define ${\cal H}_{\star}^*(E, F)^H$ to be the equalizer of  these two arrows.

 \bigskip
\noindent {\bf Theorem 2.4:}  {\sl Let $H$ be a Hopf algebra. Let $E$  be a Hopf algebra in the category $\YD_H^H$ and let $F$  be a Radford $E$-comodule algebra. The two maps 
$\iota : F \ot H \lr (F \ot H) \cross (F \ot E)$ and \hbox{$\pi : (F \ot H) \cross (F \ot E) \lr  (F \ot E)$} defined by
$$\iota (X) = (X, 1) \quad \ {\sl and }  \quad \pi (X, Y) = Y$$ give rise to the following exact
sequence
of pointed sets, in which $\iota$ is injective:
$$ \xymatrix{{\cal H}^1{\bf C}_E(H, F)  \ar@<0ex>[r]^{ \iota  \ \ }
 & 
\ {\cal H}^1(H{\star} E, F)  \ \  \ar@<0ex>[r]^{\ \ \pi  } 
& \ {\cal H}^1_{\star}(E, F)^H.}$$ 
}

\Dem We first show that $\iota $ induces an injective map, still denoted by $\iota $, from ${\cal H}^1{\bf C}_E(H, F)$
to  $F^{\lcross}\big\backslash (\ZC^1(H, F) \sqtimes \ZC_{\star}^1(E, F))$, which is isomorphic to
${\cal H}^1(H{\star} E, F)$ by Theorem 2.2. To this end take $X = x \ot h$  an element of $\ZC^1{\bf C}_E(H, F)$. 
Notice that the $E$-coinvariance of $X$ is equivalent
to the relation $\varrho_{F,E}(x) \ot h = x\ot 1 \ot h$, which is itself equivalent to the Compatibility relation  2.3 between $X$ and $1 \ot 1$. So $[X, 1]$ belongs to  $\ZC^1(H, F) \sqtimes \ZC_{\star}^1(E, F)$.

Let $X$ and $X'$ be two elements of $\ZC^1{\bf C}_E(H, F)$ such that $X' = X \leftharpoonup x$ with 
$x \in (F^E)^{\lcross}$. One has $[X, 1] \leftharpoonup x = [X \leftharpoonup x, 1 \leftharpoonup x]$.
But since  $x$ is $E$-coinvariant, one also has $1 \leftharpoonup x = (x^{-1}\ot 1)\varrho_{F,E}(x) = 1$. As a consequence,  $\iota $ 
induces a map on the quotients.

Suppose now that $X$ and $X'$ have the same image by $\iota$. This means that there exists an element 
$x \in F^{\lcross}$ such that $[X, 1]\leftharpoonup x = [X',1]$, which is equivalent to the two equalities
$X \leftharpoonup x = X'$ and $1 = 1 \leftharpoonup x$. The latter equality reflects that $x$ belongs to 
the group $(F^E)^{\lcross}$. Therefore $X$ and $X'$ are cohomologous in $\ZC^1{\bf C}_E(H, F)$. Hence $\iota$ is injective.

\medskip
We show that the image of $\pi$ lies in ${\cal H}^1_{\star}(E, F)^H$. To this end, let $[X,Y] \in \ZC^1(H, F) \sqtimes \ZC_{\star}^1(E, F)$
 represent an element of ${\cal H}^1(H{\star} E, F)$. Figure 2.11 proves that ``multiplying"  the right-hand side of the Compatibility relation 2.3 between $X$ and $Y$  by $X^{-1}$ leads to  $\varrho_{F {\star} E, H}(Y)$. Indeed, the first equality is a consequence of the associativity of the algebra $F$. The second 
one comes from Relation 0.4.a, whereas  the third one  reflects  $X^{-1}X = 1$. Finally, the last equality is given by Relation~0.4.d.

$$ \Figure[ 0.9 ]2.11->fig2-10;  $$
\medskip

In Figure 2.12 the left-hand side of the Compatibility relation 2.3 between $X$ and $Y$  is  ``multiplied" by $X^{-1}$ in the same manner as above. By inserting $\eta _E$ and $\eta _H$, one obtains the element
$(Y \ot 1) \leftharpoonup X$ in $(F {\star}E)\ot H$. 

$$ \Figure[ 0.8 ]2.12->fig2-11;  $$
So the elements $\varrho_{F {\star} E, H}(Y)$ and $Y \ot 1$ 
are cohomologous in $\ZC_{\star}^1{\bf T}_H(E, F)$, and therefore 
$\pi ([X,Y])$ belongs to the set ${\cal H}^1_{\star}(E, F)^H$.
\medskip

It is obvious that the image of $\iota$ is included in the preimage of $1$ by $\pi$. Conversely, take $[X,Y]$  in $ \ZC^1(H, F) \sqtimes \ZC_{\star}^1(E, F)$
such that its class in ${\cal H}^1(H{\star} E, F)$ has a trivial image
via $\pi$. This means that there exists an element $x$ in $F^{\lcross}$ such that
$Y = 1 \leftharpoonup x$. As $[X, 1 \leftharpoonup x]$ is equal to $[X\leftharpoonup x^{-1}, 1] \leftharpoonup x$, the class of $[X, Y]$ is also represented
by $[X\leftharpoonup x^{-1}, 1]$. By the remark formulated above, this means that $X\leftharpoonup x^{-1}$ is $E$-coinvariant. In other words, 
the class of $[X,Y]$ belongs to the image of~$\iota$.
\dm

\bigskip
\smallskip

\noindent {\bf 3. Computation of examples.}


\medskip
\noindent {\tencmbi 3.1. When $F = E$.} 
\medskip
\noindent {\bf  Proposition 3.1:} {\sl Let $H$ be a Hopf algebra
and $E$  be a Hopf algebra in the category $\YD_H^H$.
Then one has the equality
$$\HC^0(H{\star} E, E) = k^{\lcross},$$
and an isomorphism of pointed sets
$$\HC^1(H{\star} E, E) \cong  {\rm Gr}(H).$$
Here ${\rm Gr}(H)$ stands for the set of grouplike elements of $H$.

}

\medskip
\Dem 

\noindent {\sl $0$-cohomology level}. By definition, one immediately
 obtains that $\HC^0_{\star}(E, E)$
is equal to $k^{\lcross}$. Since the group  $k^{\lcross}$ is included in $\HC^0(H, E)$,  Theorem 2.2 implies the first assertion.
\medskip
\noindent {\sl $1$-cohomology level}.  We use the exact sequence given in Theorem 2.4. 
The computations of Example 1.2.3 in [4] can be adapted here, so the $1$-cohomology set $\HC^1_{\star}(E, E)$ is equal to $\{1\}$.
Therefore $\HC^1_{\star}(E, E)^H$  is trivial. On the other hand, let $x\ot h$ be an element of 
$\ZC^1{\bf C}_E(H, E)$. The $E$-coinvariance of $x\ot h$ is expressed by the equality $x_1 \ot h \ot x_2 = x \ot h \ot 1$.
Apply $\varepsilon _E$ to the first component of the tensor product in order to get
$x\ot h = 1 \ot \varepsilon _E(x)h$. Call $1 \ot h'$ the latter element. The cocycle condition fulfilled by $1 \ot h'$ implies 
 that $h'$ belongs to ${\rm Gr}(H)$. Conversely, starting with $h' \in {\rm Gr}(H)$, one checks that the element $1\ot h'$ is in 
$\ZC^1{\bf C}_E(H, E)$. 
Moreover $E^E$ is equal to the ground ring $k$, hence $\HC^1{\bf C}_E(H, E)$ is isomorphic to ${\rm Gr}(H)$. The exact sequence
of Theorem 2.4 becomes
$ \xymatrix{{\rm Gr}(H)  \ar@<0ex>[r]^{ \iota  \ \ \ \ \ \ }
 & 
\ {\cal H}^1(H{\star} E, F)  \ \  \ar@<0ex>[r]^{} 
& \ 1}$. The injectivity of $\iota$ allows us to conclude.
 \dm

\medskip
In the case when $E$ is equal to the ground ring $k$, one recovers the computation of $\HC^*(H, k)$ performed in Example 1.2.2 of [4].

\bigskip
\noindent {\tencmbxsl 3.2. Taft algebras.} 
Let here $k$ be a field. Fix $n$ an integer with $n \geq 2$ and $\zeta$ a primitive $n$-th root of 
unity in $k$. Consider the Taft algebra
$H_{n^2}$ generated by two elements $g$ and $h$ submitted to the relations
$g^n = 1$, $h^n = 0$ and  $ hg= \zeta gh$ ([10]).
On the generators the comultiplication, the antipode, and the counit of $H_{n^2}$ are given by
$$\matrix{\Delta (g) = g \ot g, \hfill &  \quad  \Delta (h) = h \ot g + 1\ot h, \hfill \cr
\sigma (g) = g^{n-1}, \hfill &  \quad \sigma (h) = -\zeta^{-1}g^{n-1}h, \hfill \cr
 \varepsilon(g) = 1, \hfill &  \quad \varepsilon(h) = 0. \hfill\cr}$$
Denote by $E_{n}$ the subalgebra of  $H_{n^2}$ generated by $h$. Via $\Delta_{H_{n^2}}$, the algebra $E_n$ 
is naturally endowed with a  structure of $H_{n^2}$-comodule algebra.

The Taft algebra $H_{n^2}$ is in fact a Radford product. Indeed, let $G$ be the cyclic group of order $n$ generated by $u$. 
As usually, denote by $k[G]$ the  Hopf algebra of the group $G$. Notice that 
the co-opposite algebra of $k[G]$ is a Hopf algebra, since 
$k[G]$ is cocommutative. Let $E$ be the algebra generated by $y$ submitted to the relation $y^n = 0$.
Endow $E$ with the $k[G]$-module and $k[G]$-comodule structures defined for any nonnegative integers  $i$ and $j$  by:
$$y^i{\cdot} u^j = \zeta ^{ij}y^i \quad {\rm and} \quad \varrho_{E, k[G]}(y^i) = y^i \ot u^i.$$
Then $E$ is an object of  $\YD_{k[G]}^{k[G]}$. Equip $E$ with
the maps $\Delta _E $, $\varepsilon_E$, and 
$\sigma_E$ given on $y^i \in E$ by
$$\Delta _E (y^i) = \sum _{s = 0}^i{\pmatrix{i \cr s \cr}\!}_{\! \zeta} \ y^s \ot y^{i-s}{\rm , } \quad \ \varepsilon_E(y^i) = 0, \ {\rm when} \  i > 0, \quad {\rm and} \quad\sigma_E(y ^i) = (-1)^i\zeta ^{i(i-1)\over 2}y^i.$$
The symbol ${\pmatrix{i \cr s \cr}\!}_{\! \zeta} $ denotes the $\zeta$-analogue of the binomial coefficient.
Provided with all these data, $E$~becomes a Hopf algebra in the braided category  $\YD_{k[G]}^{k[G]}$. 

One may now form the Radford product $k[G]{\star} E$, and convince oneself that it is isomorphic to the Taft algebra $H_{n^2}$ by sending 
$1{\star}y$ to $h$ and $u {\star} 1$ to $g$. In this correspondence $E$ is isomorphic to $E_n$ as an algebra.

\medskip
As an immediate consequence of Proposition 3.1, we obtain the following result, which proves
the conjecture formulated in [4]. 
It generalizes Proposition 1.6 of [4]. \medskip

\noindent {\bf  Corollary 3.2:} {\sl There is an   equality of groups
$$\HC^0(H_{n^2}, E_n) = k^{\lcross}$$ and an isomorphism of pointed sets
$$ \HC^1(H_{n^2}, E_n) \cong  \{1\ot g^k \  \vert \   0  \leq k \leq n-1\}.$$

}

\medskip

The class of Hopf algebras of the type $k[G]{\star}E$ is very large. We refer to [1] for general results and
various examples, as the quantum plane or the Lie colour algebras. 
By Proposition~3.1, one always gets  the equality
$\HC^0(k[G]{\star} E, E) = k^{\lcross}$ and 
the isomorphism $\HC^*(k[G]{\star} E, E) \cong G$. 

\bigskip
\medskip
\noindent {\tencmbxsl 3.3. Dual of group algebras.} Let $D$ be a finite group  with neutral element $1$. Denote by $k^{D}$ the $k$-free Hopf algebra over the $k$-basis
$\{ \delta _d\} _{g \in D}$, with the following structure maps:
the multiplication is given by ${\displaystyle \delta _d {\cdot} \delta _{d'} = \partial _{d,d'} \delta _{d}}$,
where $\partial _{d,d'}$ stands for the 
Kronecker symbol of $d$ and $d'$; 
the unit in~$k^{D}$ is the element 
${\displaystyle 1 = \sum _{d \in G}\delta _d}$;
the  comultiplication $\Del _{k^{D}}$ is
defined by 
${\displaystyle \Del _{k^{D}} (\delta _d) =
\sum _{uv = d}\delta _u \ot \delta _v}$;
the counit~$\varepsilon _{k^{D}}$ is
defined by $\varepsilon _{k^{D}} (\delta _d) = \partial _{d,1}1$; 
the antipode $\sigma_{k^{D}}$ sends $\delta _d$ on $\delta _{d^{-1}}$.
When $k$ is a field, then $k^{D}$ is the dual of the usual group algebra $k[{D}]$.

\medskip
In [4], we related the non-abelian Hopf cohomology theory of $k^{D}$ with the non-abelian cohomology theory of the group $D$ (Theorem 1.5). Recall briefly the definitions of the non-abelian cohomology theory $\Hr^i(D, C) $ (with $i = 0,1$) of a
 group $D$ with coefficients in  a
 $D$-group $C$.
The 0-cohomology object
$\Hr^0(D,C) $ is the group $C^{D}$ of 
invariant elements of $C$ under the action of $D$. 
The set $\Zr^1(D, C)$ of $1$-cocycles is given by
$$\Zr^1(D, C) = \{ \beta : D \lr C \ \vert \ \ \beta (dd') =
\beta (d){{}^{d}\! \bigl(}\beta (d')\bigr),
 \ \ \forall \   d, d' \in D \}.$$
The group
$C$ acts on the right on $\Zr^1(D, C)$  by $$(\beta \leftharpoonup x)(d) = x^{-1}\beta (d) \ {}^{d}\! x,$$
where $x \in C$, $\beta  \in \Zr^1(D, C)$, and $d \in G$.

The non-abelian $1$-cohomology set  $\Hr^1(D, C)$ is the left quotient
$C \big\backslash \Zr^1(D, C)$.
It is pointed with distinguished point the class of the constant map $1: D \lr C$. Notice that the non-abelian cohomology theory of groups may also be interpreted as the non-abelian cohomology theory associated to a pre-cosimplicial diagram of groups [4].
\medskip
\goodbreak
Let $F$ be a $k^{D}$-comodule algebra. Recall [4] that, for any $x\in F$, the formula
$$\varrho_{E,k^{D}}(x)  =  \sum _{d \in D} {^d\!x} \ot \delta _d$$
defines a left action of the  group $D$ on the algebra $F$, hence on the group $F^{\lcross}$.
Theorem 1.5 of [4] states an isomorphism between
$\HC^*(k^D, F)$ and $\Hr^*(D, F^{\lcross})$
\medskip

Consider now the case where $D$ is a semi-direct product of groups. Let $G$ and $A$ be two finite groups such that $G$ acts on the left on $A$. In our conventions the semi-direct group $G \sd A$ 
 is the set $G \cross A$ endowed with the
product given for $(g,a), (h, b) \in G \cross A$ by 
$$(g,a)(h, b) = (gh, {}^{h^{-1}}\!\!ab).$$

The co-opposite algebra of ${k^G}$ is still a Hopf algebra, so we can consider the braided monoidal category $\YD_{k^G}^{k^G}$. Equip $k^A$
with the $k^G$-module and $k^G$-comodule structures given by
$$ \delta _a {\cdot} \delta _{g} = \partial _{g,1} \delta _{a} \quad {\rm and} \quad  \varrho_{k^A, k^G} (\delta _a) = \sum _{h \in G}\delta _{{}^h\!a \ot \delta _h},$$
for any $a \in A$ and $g \in G$. Then $k^A$ becomes a Hopf algebra in the Yetter-Drinfeld category $\YD_{k^G}^{k^G}$ and one has an isomorphism 
of Hopf algebras
$$k^G {\star} k^A \cong k^{G \lsd A}$$ in which the element $\delta _g {\star} \delta _a$ corresponds with $\delta _{(g,a)}$.
When the action of $G$ on $A$ is trivial, then $k^A$ is equipped with the trivial $k^G$-module and $k^G$-comodule structures and the above isomorphism reduces
to $$k^G {\otimes}k^A \cong k^{G \lcross A}.$$
Theorem 1.5 of [4] and the previous Theorem 2.2 imply the following result:

\bigskip
\noindent {\bf  Proposition 3.3:} {\sl Let $G$ and $A$ be two finite groups such that $G$ acts on the left on $A$.  Let $F$ be a $k^{G \lsd A}$-comodule algebra.
There is an  equality of groups
$$\Hr^0(G \sd A, F^{\lcross}) = \Hr^0(G, F^{\lcross}) \cap \Hr^0(A, F^{\lcross}),$$
and an isomorphism of pointed sets
$$\Hr^1(G \sd A, F^{\lcross}) \cong F^{\lcross} \big\backslash (\ZC^1(k^G, F) \sqtimes \ZC_{\star}^1(k^A, F)).$$
}

\medskip
\goodbreak

\noindent {\bf 4.  Digression: The non-abelian cohomology theory for semi-direct products of groups.} With the help of
the Hopf algebra machinery, we are  able  by Proposition 3.3 to compute the non-abelian cohomology theory of a semi-direct group with 
coefficients in a group of a particular type: the group of invertible elements of a certain algebra. One should now like to provide such a calculation in the case of 
an arbitrary coefficient group. In the present section,
we show how to to do this directly, using only the basic definitions of the non-abelian cohomology theory ([6], [7]). The results we obtain are the same as those in Proposition 3.3. This closeness may be explained by a metamathematical  
fantasy suggested by several mathematicians,  the "field of one element" $\Fb_1$ (see [8]). Indeed, in this philosophy,
a group is nothing but a Hopf algebra over the field  $\Fb_1$.

\medskip

Let now $G$ and $A$ be two finite groups such that $G$ acts on the left on $A$. Let $C$ be a group. It is easy to see that $C$  is a  $G \sd A$-group if and only if
$C$ is both a $G$-group and an $A$-group such that the two actions
satisfy the following compatibility law:
$$^{{}^g}\!(^ax) = \, ^{{}^{^g\!a}}\!(^gx),$$
for any $g \in G$, $a \in A$, and $x \in C$. This is to say that the $G$-group $A$ acts on the left on $C$ or, in other words, that the canonical map $A\cross C \lr C$
is a morphism of $G$-sets. The dictionary between the several actions is given by
$^{(g,a)}\!(x) = \, ^{{}^g}\!(^ax)$, for $(g,a) \in  G \sd A$ and  $x \in C$.

We define now the set $\Zr^1(G, C) \sqtimes \Zr^1(A, C)$. It
elements $[\gamma, \alpha]$ are couples  in $\Zr^1(G, C) \cross \Zr^1(A, C)$ such that $\gamma$ and $\alpha$ satisfy the compatibility relation 
$$\gamma(g) \ ^g\!\alpha(a) = \alpha(^g\!a)\ ^{{^g\!a}}\gamma(g).$$
Observe that $\Zr^1(G, C) \sqtimes \Zr^1(A, C)$ is a pointed set with $[1, 1]$ as distinguished point.
The Reader easily checks that the group $C$ acts on the right  on  $\Zr^1(G, C) \sqtimes \Zr^1(A, C)$ by
$$[\gamma, \alpha] \leftharpoonup c = [\gamma \leftharpoonup c,\alpha \leftharpoonup c],$$
for $[\gamma, \alpha] \in  \Zr^1(G, C) \sqtimes \Zr^1(A, C)$ and $c \in C$.

\medskip
The group counterpart of Theorem 2.2 is the following result:

\medskip
\noindent {\bf  Proposition 4.1:} {\sl Let $G$ and $A$ be two finite groups such that $G$ acts on the left on $A$.  Let $C$ be a $G \sd A$-group.
There is an  equality of groups
$$\Hr^0(G \sd A, C) = \Hr^0(G, C) \cap \Hr^0(A, C),$$
and an isomorphism of pointed sets
$$\Hr^1(G \sd A, C) \cong C \big\backslash (\Zr^1(G, C) \sqtimes \Zr^1(A, C)).$$
}

\Dem The proof in the $0$-level is a direct consequence of the definitions. On the first level, it is enough to exhibit an explicit isomorphism
between $\Zr^1(G \sd A, C)$ and $\Zr^1(G, C) \sqtimes \Zr^1(A, C)$. The correspondence 
$\beta \longmapsto [\gamma, \alpha]$ defined by $\gamma (g) = \beta (g, 1)$ and $ \alpha (a) = \beta (1,a)$,
where $g \in G$ and $a \in A$,
realize this isomorphism. Conversely, a compatible couple $[\gamma, \alpha]$ gives raise to the element $\beta \in \Zr^1(G \sd A, C)$
defined by $\beta (g,a) = \gamma (g) \ ^g\!\alpha (a),$ for any $(g,a) \in G \sd A$. \dm

\medskip
\goodbreak
Notice that  $k^C$ is naturally equipped with a $k^G{\star}k^A$-comodule algebra
structure when  $C$ is a $G \sd A$-group. 
Via Proposition 3.3 and [4, Theorem 1.5], 
one obtains an  equality of groups
$$\Hr^0(G \sd A, (k^C)^{\lcross}) = \Hr^0(G, (k^C)^{\lcross}) \cap \Hr^0(A, (k^C)^{\lcross}),$$
and an isomorphism of pointed sets
$$\Hr^1(G \sd A, (k^C)^{\lcross}) \cong (k^C)^{\lcross} \Bigl\backslash (\Zr^1(G, (k^C)^{\lcross}) \sqtimes \Zr^1(A, (k^C)^{\lcross})).$$
Here we could replace $\ZC_{\star}^1(k^A, k^C)$ by $\ZC^1(k^A, k^C)$, using the fact that the braiding in $\YD^{k^G}_{k^G}$ is the usual flip
for the objects of type $k^D$, for any $G$-group $D$.
This latter result differs from that stated in Proposition 4.1, since the group of invertible elements of $k^C$ is in general not isomorphic to $C$ (unless $k = \Fb_1$!).

\medskip
We end this digression by observing that the group counterpart of
 Theorem 2.4 is the exact sequence of pointed sets 
 $$\xymatrix{   \Hr^1(G, C^A)
\ar@<0ex>[r]^{\iota \ \ }
&  \Hr^1(G \sd A, C)
\ar@<0ex>[r]^{\pi}
&  \Hr^1(A , C)^G}$$ (with $\iota$
 injective). Through the isomorphism described in  Proposition 4.1, the maps $\iota$ and $\pi$ may be simply expressed by the formulae
$\iota (\gamma) = [\gamma, 1]$ and $\pi  ([\gamma, \alpha]) = \alpha$.
This cohomological sequence is by Serre's theory [7] recovered as the one associated with the
 exact sequence of groups $1 \lr A \lr G \sd A \lr G \lr 1$.

\bigskip
\bigskip
\bigskip

\noindent {R{\eightrom EFERENCES}}

\medskip
\medskip

\item{[1]} D. F{\eightrom ISCHMAN}, S.  M{\eightrom ONTGOMERY},
A Schur double centralizer theorem for cotriangular Hopf algebras and generalized Lie algebras, {\it J. Algebra}
{\bf 168}, ({\oldstyle 1994}), n$^{\hbox{{\fiverom o}}}$ $\!$2, 594 -- 614.
\medskip


\item{[2]} S. M{\eightrom AJID}, 
Cross products by braided groups and bosonization, 
{\it J. Algebra}
{\bf 163}, ({\oldstyle 1994}), n$^{\hbox{{\fiverom o}}}$ $\!$1, 165 -- 190.
\medskip

\item{[3]} Ph. N{\eightrom USS}, M.  W{\eightrom AMBST}, 
Non-Abelian Hopf Cohomology, {\it J. Algebra}
{\bf 312}, ({\oldstyle 2007}), n$^{\hbox{{\fiverom o}}}$ $\!$2, 733 -- 754.
\medskip

\item{[4]} Ph. N{\eightrom USS}, M.  W{\eightrom AMBST},
Non-abelian Hopf cohomology II - The general case, {\it J. Algebra}
{\bf 319}, ({\oldstyle 2008}), n$^{\hbox{{\fiverom o}}}$ $\!$11, 4621 -- 4645.
\medskip

\item{[5]} D. E. R{\eightrom ADFORD},
The structure of Hopf algebras with a projection, {\it J. Algebra}
{\bf 92}, ({\oldstyle 1985}), n$^{\hbox{{\fiverom o}}}$~2, 322 -- 347.
\medskip

\item{[6]} J.-P. S{\eightrom ERRE}, 
{\it  Corps locaux}, Troisi\`eme \'edition corrig\'ee,
Hermann, Paris  ({\oldstyle 1968}).

\medskip

\item{[7]} J.-P. S{\eightrom ERRE}, 
{\it  Galois cohomology}, Springer-Verlag, Berlin --
Heidelberg  ({\oldstyle 1997}). Translated from
{\it  Cohomologie galoisienne},
Lecture Notes in Mathematics 5, Springer-Verlag, Berlin --
Heidelberg -- New York ({\oldstyle 1973}).

\medskip

\item{[8]}  Chr. S{\eightrom OUL\'E}, 
Les vari\'et\'es sur le corps \`a un \'el\'ement, {\it Mosc. Math. J.} 4 (2004), n$^{\hbox{{\fiverom o}}}$ $\!$1, 217--244, 312. 
\medskip

\item{[9]} M. E. S{\eightrom WEEDLER},  
 Cohomology of algebras over Hopf algebras, {\it Trans. Amer. Math. Soc.} {\bf 133}  ({\oldstyle 1968}), 205~--~239.

\medskip
\item{[10]} E. J. T{\eightrom AFT},
The order of the antipode of finite-dimensional Hopf algebra,
{\it Proc. Nat. Acad. Sci. USA} {\bf 68} ({\oldstyle 1971}), n$^{\hbox{{\fiverom o}}}$ $\!$11, 2631  -- 2633.
\medskip

\item{[11]} D. N. Y{\eightrom ETTER},
Quantum groups and representations of monoidal categories,
{\it Math. Proc. Camb. Phil. Soc.} {\bf 108} ({\oldstyle 1990}), n$^{\hbox{{\fiverom o}}}$ $\!$2, 261  -- 290.
\medskip

\bye